\documentclass[12pt]{amsart}
\usepackage{amsmath,amssymb,amsthm,amsfonts}
\usepackage{url}
\usepackage[margin=1in]{geometry}
\usepackage{cleveref}

\usepackage{setspace}
\theoremstyle{plain}
\newtheorem{theorem}{Theorem}[section]
\newtheorem{proposition}[theorem]{Proposition}
\newtheorem{lemma}[theorem]{Lemma}
\newtheorem{corollary}[theorem]{Corollary}
\theoremstyle{definition}
\newtheorem{definition}[theorem]{Definition}
\theoremstyle{remark}
\newtheorem{remark}[theorem]{Remark}

\newcommand{\Hn}{\mathbb{H}^n}
\newcommand{\Cn}{\mathbb{C}^n}
\newcommand{\R}{\mathbb{R}}
\newcommand{\N}{\mathbb{N}}
\newcommand{\Rstar}{\mathbb{R}^*}

\title[Characterisations of the Strichartz Fourier Transform]
{Some Characterisations
of the Strichartz Fourier Transform on the Heisenberg Group}

\author{R. Lakshmi Lavanya}
\address{Department of Mathematics, Indian Institute of Science Education and Research (IISER) Tirupati, Yerpedu, Tirupati 517619, Andhra Pradesh, India.}
\email{rlakshmilavanya@labs.iisertirupati.ac.in}

\date{}

\begin{document}

\maketitle

\begin{abstract}
Recently, Thangavelu introduced the Strichartz Fourier transform on the Heisenberg group. Around the same time, Krukowski characterised the classical Fourier transform on various settings, assuming a priori that the map is an integral transform, based on structural properties and a normalisation. In this paper, we obtain analogous results for the Strichartz Fourier transform, without assuming the characterised map is an integral transform. The spectral property is formulated jointly for the sublaplacian and central vector field, and the translation property uses a twisted form of left translation. The joint spectral property, the twisted translation property and a normalisation at one explicit test function, together characterise the transform on each ray of the Heisenberg fan. The twisted translation property, paired with a Dirac delta property, gives another characterisation. We characterise the radial kernel on each ray by a Laguerre-type differential equation. This gives a characterisation of the Hankel transform.
\end{abstract}

\noindent{\small\textbf{2020 Mathematics Subject Classification.} Primary 43A85, 42C05; Secondary 33C45, 43A80.}

\section{Introduction}

Characterising the Fourier transform, and its analogues on other groups, by a small number of structural properties rather than by an explicit formula is a recurring theme in harmonic analysis. In 2008, Alesker, Artstein-Avidan, and Milman \cite{AAM} showed that the Fourier transform is, up to a linear change of variables and complex conjugation, the unique transformation on tempered distributions on $\R^n$ exchanging convolution and pointwise products. Jaming \cite{Jaming} later characterised the Fourier transform using only the weaker property of taking convolutions to pointwise products. Kumar and Sivananthan (\cite{KS}) obtained an analogous result for arbitrary compact groups, and Lavanya and Thangavelu (\cite{LT1,LT2}) did so for the operator-valued group Fourier transform on the Heisenberg group $\Hn$. Recently, Achar, Biswas and Manna obtained characterisations of the Euclidean Fourier transform (\cite{ABM1}), and the group Fourier transform on the Heisenberg group (\cite{ABM2}), both for vector-valued functions. 

In 2024, Krukowski \cite{Kru} went beyond the convolution property, obtaining several further characterisations of the Fourier transform. On the real line, these are based instead on a time differentiation property. On locally compact abelian groups, where differentiation is no longer available, a time shift property plays the analogous role. A further result, of a different kind, characterises the Hankel transform by a Bessel differential property. Each of these results is paired with a suitable normalisation hypothesis. He also revisited the convolution property itself, giving a shorter, representation-theoretic proof that it characterises the Fourier transform on compact groups, simplifying the earlier argument of Kumar and Sivananthan (\cite{KS}).

Meanwhile, a new scalar-valued Fourier transform was being developed on the Heisenberg group. Thangavelu \cite{Th23} introduced the \emph{Strichartz Fourier transform}, arising from the joint spectral theory of the sublaplacian $\mathcal L$ and the central invariant vector field $T=\partial/\partial t$. It plays a role for $\Hn$ analogous to that of the Helgason Fourier transform on rank one Riemannian symmetric spaces of noncompact type. Unlike the classical, operator-valued group Fourier transform on $\Hn$, the Strichartz Fourier transform is complex-valued. 

In this paper, we obtain characterisations of the Strichartz Fourier transform in the spirit of \cite{Kru}, built from properties intrinsic to $\Hn$ rather than by adapting the proofs of \cite{Kru} directly. The scalar nature of the Strichartz Fourier transform simplifies matters considerably compared to the operator-valued group Fourier transform on $\Hn$ \cite{Khn}. In the current context, the Strichartz Fourier transform of a function is already scalar-valued, so no operator-valued kernel need be extracted from a bare continuous linear map. The relevant structure instead is the joint spectral theory of the sublaplacian and the central invariant vector field, together with a twisted translation property.

We obtain three main results. The first shows that if a bounded linear map $B:L^1(\Hn)\to C^b(\Omega\times\Cn)$, where $\Omega$ denotes the Heisenberg fan, satisfies the joint spectral property of $(\mathcal L,T)$ and the twisted translation property, and agrees with the Strichartz Fourier transform at one explicit test function for each spectral parameter, then it coincides with it on each ray of the fan. The second shows that the twisted translation property alone, combined instead with convergence to the Strichartz kernel along an approximate identity, already determines $B$, with no spectral hypothesis needed. Neither result is a special case of the other - the first exchanges a stronger, but more localised, normalisation for the additional joint spectral hypothesis. The final result gives a ray-by-ray characterisation of the radial kernel of the Strichartz Fourier transform, from a Laguerre-type ordinary differential equation, and a normalisation.

The paper is organised as follows. \Cref{sec:kru} recalls the results of \cite{Kru} that motivate the theorems below. \Cref{sec:prelims} recalls the preliminaries on the Strichartz Fourier transform from \cite{Th23}. \Cref{sec:necessary} introduces the structural properties used in the characterisations - the joint spectral property, the twisted translation property, the Dirac delta property, the Laguerre differential property, and  verifies that the Strichartz Fourier transform satisfies each of them. \Cref{sec:spectral} combines the joint spectral property with the twisted translation property and a pointwise normalisation. \Cref{sec:translation} relies on the twisted translation property alone, paired with the Dirac delta property. \Cref{sec:laguerre} characterises the radial kernel of the  Strichartz Fourier transform on each ray in terms of the Laguerre differential equation. This result gives a characterisation of the Hankel transform on the limiting ray.

\section{Characterisations of the Fourier transform by structural properties}
\label{sec:kru}

We first recall, for comparison, the three results of \cite{Kru} that motivate the theorems of this paper.

\begin{definition}[\cite{Kru}]\label{def:kru-timediff}
Let $C^b(\R)$ be the bounded continuous functions on $\R$, and $C_c^1(\R)$ the compactly supported $C^1$ functions. An integral transform $$B:L^1(\R)\to C^b(\R), \ B(f)(y)=\int_\R K(x,y)f(x)\,dx,$$ with $K\in C^b(\R\times\R)$ differentiable in the first variable, satisfies the \emph{time differentiation property} if
$$
B(f')(y)=2\pi i y\,B(f)(y),\qquad \forall\,f\in C_c^1(\R),\ y\in\R.
$$
\end{definition}

\begin{definition}[\cite{Kru}]\label{def:kru-deltaseq}
A sequence $(\delta_n)\subset C_c^1(\R)$ is a \emph{Dirac delta sequence} if

(i) $\operatorname{supp}\delta_n\subset[-1/n,1/n]$,

(ii) $\delta_n$ is nonnegative and symmetric, i.e.\ $\delta_n(x)=\delta_n(-x)$,

(iii) $\int_\R\delta_n(x)\,dx=1$.
\end{definition}

\begin{definition}[\cite{Kru}]\label{def:kru-deltaprop-R}
$B$ as in \Cref{def:kru-timediff} satisfies the \emph{Dirac delta property} if
$$
B(\delta_n)(y)\to1,\qquad \forall\,y\in\R,
$$
for some Dirac delta sequence $(\delta_n)$.
\end{definition}

\begin{theorem}[\cite{Kru}]\label{thm:kru1}
Let $B:L^1(\R)\to C^b(\R)$ be an integral transform as in \Cref{def:kru-timediff}. If $B$ satisfies the time differentiation property and the Dirac delta property, then $B=\mathcal F$, where $\mathcal F(f)(y):=\int_\R f(x)e^{-2\pi ixy}\,dx$ is the Fourier transform on $\R$.
\end{theorem}

For a locally compact abelian group $G$ with dual $\widehat G$, let $N_0$ denote the set of open neighbourhoods of the identity of $G$, and $L_xf(t):=f(t-x)$ the translate of $f$ by $x\in G$.

\begin{definition}[\cite{Kru}]\label{def:kru-timeshift}
A linear bounded $B:L^1(G)\to C^b(\widehat G)$ satisfies the \emph{time shift property} if
$$
B(L_xf)(\chi)=\overline{\chi(x)}\,B(f)(\chi),\qquad \forall\,f\in L^1(G),\ x\in G,\ \chi\in\widehat G.
$$
\end{definition}

\begin{definition}[\cite{Kru}]\label{def:kru-deltafam}
A family $(\delta_U)_{U\in N_0}\subset C_c(G)$ is a \emph{Dirac delta family} if, for every $U\in N_0$,

(i) $\delta_U$ is continuous and compactly supported,

(ii) $\delta_U$ is nonnegative and symmetric, i.e.\ $\delta_U(x)=\delta_U(-x)$,

(iii) $\int_G\delta_U(x)\,dx=1$.
\end{definition}

\begin{definition}[\cite{Kru}]\label{def:kru-deltaprop-G}
$B$ as in \Cref{def:kru-timeshift} satisfies the \emph{Dirac delta property} if
$$
B(\delta_U)(\chi)\to1,\qquad \forall\,\chi\in\widehat G,
$$
for some Dirac delta family $(\delta_U)$.
\end{definition}

\begin{theorem}[\cite{Kru}]\label{thm:kru3}
If a linear bounded $B:L^1(G)\to C^b(\widehat G)$ satisfies the time shift property and the Dirac delta property, then $B=\mathcal F$, the Fourier transform on $G$.
\end{theorem}

Let $J_\alpha$ be the Bessel function of the first kind of order $\alpha$, and let $\mathcal H_\alpha:C_c(\R_+)\to C_0(\R_+)$, $\mathcal H_\alpha(f)(y):=\int_{\R_+}J_\alpha(yx)f(x)x\,dx$, be the Hankel transform of order $\alpha$.

\begin{definition}[\cite{Kru}]\label{def:kru-bessel}
An integral transform $B(f)(y)=\int_{\R_+}K(yx)f(x)x\,dx$, with $K\in C^2(\R_+)$ bounded, satisfies the \emph{Bessel differential property} of order $\alpha>0$ if
$$
B\Big(f''+\tfrac1x f'-\tfrac{\alpha^2}{x^2}f\Big)(y)=-y^2\,B(f)(y),\qquad \forall\,f\in C_c^2(\R_+).
$$
\end{definition}

\begin{theorem}[\cite{Kru}]\label{thm:kru5}
Let $B$ be as in \Cref{def:kru-bessel}. If $B$ satisfies the Bessel differential property of order $\alpha>0$, and $B(f_*)(y_*)=\mathcal H_\alpha(f_*)(y_*)\ne0$ for some $f_*\in C_c^2(\R_+)$ and $y_*\in\R_+$, then $B=\mathcal H_\alpha$.
\end{theorem}

\Cref{thm:kru1,thm:kru3,thm:kru5} are the model, respectively, for \Cref{thm:spectral,thm:translation,thm:laguerre} below. As explained in the introduction, we do not reproduce the proofs of \cite{Kru}, since they rely on features specific to $\R$, to locally compact abelian groups, and to the Bessel equation. Only the shape of these results is used to guide the statements that follow.

In \cite{Kru}, the unknown transform under consideration is denoted by $T$. We reserve $T$ throughout this paper for the central invariant vector field on $\Hn$ (\Cref{sec:prelims}), and write $B$ instead, departing here from the notation of \cite{Kru}. We write $\mathcal F$ for the Fourier transform on the relevant space in each result above. The notation $\widehat{\ \cdot\ }$ is reserved throughout the rest of this paper for the Strichartz Fourier transform of \Cref{sec:prelims}.

\section{Preliminaries}
\label{sec:prelims}

In this section, we recall the definitions we need from \cite{Th23}, to which we refer for further details.

Let $\Hn=\Cn\times\R$ be the $(2n+1)$-dimensional Heisenberg group, with group law
$$
(z,t)(w,s)=\Big(z+w,\ t+s+\tfrac12\operatorname{Im}(z\cdot\bar w)\Big),\qquad (z,t),(w,s)\in\Hn.
$$
This is a step two nilpotent, unimodular Lie group, with Haar measure equal to Lebesgue measure $dz\,dt$ on $\Cn\times\R$. For $(w,s)\in\Hn$, write $L_{(w,s)}f(z,t):=f((w,s)^{-1}(z,t))$ for the left translate of $f$ by $(w,s)$, and
\[
(f*g)(z,t):=\int_{\Hn}f(w,s)\,g((w,s)^{-1}(z,t))\,dw\,ds
\]
for the convolution of two functions on $\Hn$.

The Heisenberg Lie algebra has left-invariant vector fields
$$
X_j=\frac{\partial}{\partial x_j}+\frac12 y_j\frac{\partial}{\partial t},\qquad
Y_j=\frac{\partial}{\partial y_j}-\frac12 x_j\frac{\partial}{\partial t},\qquad
T=\frac{\partial}{\partial t},\qquad j=1,\dots,n,
$$
and the sublaplacian is $\mathcal L=-\sum_{j=1}^n(X_j^2+Y_j^2)$. The corresponding right-invariant vector fields are
$$
\widetilde X_j=\frac{\partial}{\partial x_j}-\frac12 y_j\frac{\partial}{\partial t},\qquad
\widetilde Y_j=\frac{\partial}{\partial y_j}+\frac12 x_j\frac{\partial}{\partial t},\qquad
\widetilde T=T,\qquad j=1,\dots,n.
$$
That is, $(\widetilde X g)(w) = \frac{d}{dr}\Big|_{r=0} g(e^{rX}w)$, where $e^{rX}$ is the one-parameter subgroup generated by $X$. For instance, $e^{rX_j} = (re_j,0)$, so that $(re_j,0)(z,t) = (z+re_j,\,t-\tfrac12 ry_j)$.

Let $\widetilde{\mathcal L} := -\sum_{j=1}^n(\widetilde X_j^2+\widetilde Y_j^2)$ denote the right-invariant sublaplacian. Let $\Delta_z$ denote the Laplacian on $\mathbb C^n\cong\mathbb R^{2n}$, and let\\ $N := \sum_{j=1}^n(x_j\partial_{y_j}-y_j\partial_{x_j})$ denote the rotation operator. A direct computation then gives
\[
\mathcal L = -\Delta_z + N\partial_t - \tfrac14|z|^2\partial_t^2,\qquad \widetilde{\mathcal L} = -\Delta_z - N\partial_t - \tfrac14|z|^2\partial_t^2.
\]
Hence $\mathcal L-\widetilde{\mathcal L}=2N\partial_t$, and $\widetilde{\mathcal L}g=\mathcal Lg$ whenever $g(z,t)$ is radial in $z$, since $N$ annihilates functions radial in $z$.

Let $\mathcal S(\Hn)$ denote the Schwartz space on the underlying manifold $\Hn\cong\R^{2n+1}$, and $\mathcal S(\Cn)$ the Schwartz space on $\Cn\cong\R^{2n}$.

\begin{lemma}\label{lem:leftright}
Let $f\in\mathcal S(\Hn)$. Let $g\in C^\infty(\Hn)$ have all derivatives of polynomial growth, and let $X$ be one of $X_j,Y_j,T$. Then $(Xf)*g=f*(\widetilde Xg)$, and consequently $(\mathcal Lf)*g=f*(\widetilde{\mathcal L}g)$.
\end{lemma}

\begin{proof}
Since $X$ is divergence free for Haar measure and $f$ is rapidly decreasing, integration by parts gives
\begin{eqnarray*}
\big((Xf)*g\big)(x) & = & \int_{\Hn}(Xf)(y)\,g(y^{-1}x)\,dy\\
& = & -\int_{\Hn}f(y)\,X_y\big[g(y^{-1}x)\big]\,dy,
\end{eqnarray*}
and
\begin{eqnarray*}
X_y[g(y^{-1}x)] & = & \frac{d}{dr}\Big|_{r=0}g\big((ye^{rX})^{-1}x\big)\ =\ \frac{d}{dr}\Big|_{r=0}g\big(e^{-rX}y^{-1}x\big)\\
& = & -(\widetilde Xg)(y^{-1}x).
\end{eqnarray*}
Hence $((Xf)*g)(x)=(f*\widetilde Xg)(x)$. Applying this twice for each of $X_j,Y_j$ and summing gives the statement for $\mathcal L$.
\end{proof}

For $f$ on $\Hn$ and $\lambda\in\R$, write
$$
f^\lambda(z):=\int_{-\infty}^\infty e^{i\lambda t}f(z,t)\,dt
$$
for the partial Fourier transform of $f$ in the central variable. For $f,g$ on $\Cn$, the $\lambda$-twisted convolution is
$$
f*_\lambda g(z):=\int_{\Cn} f(z-w)g(w)e^{\frac{i\lambda}{2}\operatorname{Im}(z\cdot\bar w)}\,dw.
$$
Substituting $w\mapsto z-w$ and using $\operatorname{Im}(z\cdot\overline{(z-w)})=-\operatorname{Im}(z\cdot\bar w)$ gives the identity
\begin{equation}\label{eq:symmetry}
f*_\mu g=g*_{-\mu}f, \qquad \mu\in\R,
\end{equation}
and $(f*g)^\mu=f^\mu*_\mu g^\mu$ for the convolution on $\Hn$. The special Hermite operator $L_\lambda$ is defined by $\mathcal L\big(f(\cdot)e^{-i\lambda t}\big)=e^{-i\lambda t}L_\lambda f(\cdot)$, so that $(\mathcal Lh)^\lambda=L_\lambda(h^\lambda)$ for every $\lambda$. Explicitly, on $\Cn\cong\R^{2n}$,
$$
L_\lambda=-\Delta+\tfrac14\lambda^2|z|^2-i\lambda N
$$
(\cite[Chapter 1]{Th93}). Since $N$ annihilates radial functions, on a radial function $g(z)=g_0(|z|)$ we have
$$
L_\lambda g_0(r)=-g_0''(r)-\frac{2n-1}{r}g_0'(r)+\frac{\lambda^2}{4}r^2g_0(r),\qquad r=|z|.
$$
The operator $L_\lambda$ is essentially self-adjoint on $C_c^\infty(\Cn)$, and has pure point spectrum $\{(2k+n)|\lambda|:k\in\N\}$, each eigenvalue of infinite multiplicity (\cite[Chapter 1]{Th93}). 

Left invariance gives $\mathcal L(f*g)=f*\mathcal Lg$. Taking $F:=f^\mu$ and $G:=g^\mu$, and using $(f*g)^\mu=f^\mu*_\mu g^\mu$ gives,
\[
L_\mu(F*_\mu G) = F*_\mu(L_\mu G), \qquad F,G\in\mathcal S(\Cn).
\]
The standard smoothness and rapid decay properties of the Fourier transform in the central variable, together with those of differential operators with polynomial coefficients, show that the maps $f\mapsto f^\lambda$, $\mathcal L$, and $T$ each preserve $\mathcal S(\Hn)$(\cite{Rud}). Similarly, for each $\lambda\in\Rstar$, $L_\lambda$ preserves $\mathcal S(\Cn)$. Its eigenvalues $(2m+n)|\lambda|$, $m\in\N$, are spaced $2|\lambda|$ apart. For every $k\in\N$, write $E_k^{-\lambda}$ for the $k$-th eigenspace $\{h\in L^2(\Cn): L_\lambda h = (2k+n)|\lambda|h\}$. The resolvent $(L_\lambda-(2k+n)|\lambda|)^{-1}$ is then bounded on $(E_k^{-\lambda})^\perp$, and preserves $\mathcal S(\Cn)\cap(E_k^{-\lambda})^\perp$. This uses the characterisation of $\mathcal S(\Cn)$ in terms of rapid decay of the special Hermite coefficients. Now, $\Hn\cong\R^{2n+1}$ as a manifold, and the group law is polynomial in coordinates. So $\mathcal S(\Hn)$ is closed under convolution product (\cite{Th98}).

For $\delta\ge -1/2$, let $L_k^\delta$ denote the Laguerre polynomial of type $\delta$, and for $\lambda\in\Rstar=\R\setminus\{0\}$, let the Laguerre function be
$$
\varphi_{k,\lambda}^{n-1}(z):=L_k^{n-1}\big(\tfrac12|\lambda||z|^2\big)e^{-\frac14|\lambda||z|^2},\qquad z\in\Cn.
$$
By \cite[Chapter 1]{Th93}, $\varphi_{k,\lambda}^{n-1}$ is an eigenfunction of $L_\lambda$, with eigenvalue $(2k+n)|\lambda|$, which depends on $\lambda$ only through $|\lambda|$. Equivalently, on the radial variable $\rho=|z|$, $\varphi_{k,\lambda}^{n-1}$ satisfies the ordinary differential equation
\begin{equation}\label{eq:laguerreODE}
-\varphi''(\rho)-\frac{2n-1}{\rho}\varphi'(\rho)+\frac{\lambda^2}{4}\rho^2\varphi(\rho)=(2k+n)|\lambda|\,\varphi(\rho),\qquad \rho>0,
\end{equation}
obtained by the radial reduction of $L_\lambda$ above. The functions
$$
e_{k,\lambda}^{n-1}(z,t):=e^{i\lambda t}\varphi_{k,\lambda}^{n-1}(z)
$$
are joint eigenfunctions of $\mathcal L$ and $T$:
$$
\mathcal L e_{k,\lambda}^{n-1}=(2k+n)|\lambda|\,e_{k,\lambda}^{n-1},\qquad -iTe_{k,\lambda}^{n-1}=\lambda\,e_{k,\lambda}^{n-1}.
$$
The first identity holds because $\varphi_{k,\lambda}^{n-1}$ is radial, so $N$ annihilates it, and $\mathcal L$ acts on $e^{i\lambda t}\varphi$ as $-\Delta+\tfrac14\lambda^2|z|^2$.

The \emph{Heisenberg fan} $\Omega$ is the union of the rays
\[
R_k=\{(\lambda,\tau)\in\R^2:\tau=(2k+n)|\lambda|\}, \qquad k=0,1,2,\dots,
\]
together with the limiting ray $R_\infty=\{(0,\tau):\tau\ge 0\}$. For each $k$, write
\[
R_k^*:=R_k\setminus\{(0,0)\} = \{(\lambda,(2k+n)|\lambda|)\in R_k : \lambda\ne0\},
\]
and set $\Omega^*:=\bigcup_kR_k^*$. For $a=(\lambda,(2k+n)|\lambda|)\in R_k^*$, write $e_a$ for $e_{k,\lambda}^{n-1}$. The results below all concern points $a\in\Omega^*$ and the limiting ray $R_\infty$ plays no further role.

For $f\in L^1\cap L^2(\Hn)$, the \emph{Strichartz Fourier transform} of $f$ is, for $a=(\lambda,(2k+n)|\lambda|)\in R_k^*$,
\begin{eqnarray*}
\widehat f(a,z) & := & \int_{\Hn} f(w,s)\,e_a\big((w,s)^{-1}(z,0)\big)\,dw\,ds\ =\ (f*e_a)(z,0)\\
& = & f^{-\lambda}*_{-\lambda}\varphi_{k,\lambda}^{n-1}(z),\qquad z\in\Cn.
\end{eqnarray*}
On the limiting ray $R_\infty$, a separate definition is given in \cite{Th23}, involving the Bessel function $J_{n-1}$ rather than the Laguerre functions $\varphi_{k,\lambda}^{n-1}$ above. As noted already, $R_\infty$ plays no role in what follows, so we do not record it.

For each $\lambda\in\Rstar$, the $k$-th special Hermite eigenspace
$$
\mathcal E_k^\lambda:=\big\{h\in L^2(\Cn):L_{-\lambda} h=(2k+n)|\lambda|h\big\},
$$
is a closed, infinite-dimensional subspace of $L^2(\Cn)$. The spaces $\mathcal E_k^\lambda$, $k\in\N$, are mutually orthogonal, and for $g\in L^2(\Cn)$,
$$
g=(2\pi)^{-n}|\lambda|^n\sum_{k=0}^\infty g*_{-\lambda} \varphi_{k,\lambda}^{n-1}
$$
is the special Hermite expansion of $g$ at parameter $-\lambda$, the $k$-th term being the orthogonal projection $P_k(g)$ of $g$ onto $\mathcal E_k^\lambda$ \cite{Th23}, applied at parameter $-\lambda$. This sign makes $\widehat f(a,\cdot)\in\mathcal E_k^\lambda$, since
\begin{align*}
L_{-\lambda}\big(F*_{-\lambda}\varphi_{k,\lambda}^{n-1}\big) &= F*_{-\lambda}L_{-\lambda}\varphi_{k,\lambda}^{n-1}\\
&= (2k+n)|\lambda|\,F*_{-\lambda}\varphi_{k,\lambda}^{n-1}.
\end{align*} In particular, since $\varphi_{k,\lambda}^{n-1}\in\mathcal E_k^\lambda$, the $k$-th eigenspace projection fixes it, giving
\begin{equation}\label{eq:special-case}
\varphi_{k,\lambda}^{n-1}*_{-\lambda}\varphi_{k,\lambda}^{n-1}=(2\pi)^n|\lambda|^{-n}\varphi_{k,\lambda}^{n-1}.
\end{equation}
We equip $C^b(\Omega\times\Cn)$ with the supremum norm, and call a linear map \[
B:L^1(\Hn)\to C^b(\Omega\times\Cn)
\]
\textit{bounded} if
$$
\|B\|_{\mathrm{op}}:=\sup_{\|f\|_1\le 1}\ \sup_{(a,z)\in\Omega\times\Cn}|B(f)(a,z)|<\infty,
$$
and similarly for maps into $C^b(\Rstar)$.

\section{Necessary properties of the Strichartz Fourier transform}
\label{sec:necessary}

Following \cite{Kru}, we now introduce, for a linear map $B$ on $\Hn$, a small number of structural properties native to $\Hn$. Then we verify that the Strichartz Fourier transform, and its radial kernel on each ray $R_k$, satisfy all of these properties. In the subsequent sections, we prove that these properties in various combinations, are also sufficient to recover it.

\begin{definition}[\cite{Khn}]\label{def:deltaseq-Hn}
A sequence $(\delta_j) \subset C_c(\Hn)$ is a \emph{Heisenberg Dirac delta sequence} if
\begin{enumerate}
\item[(i)] $\delta_j \ge 0$,
\item[(ii)] $\displaystyle\int_{\Hn} \delta_j = 1$,
\item[(iii)] $\operatorname{supp} \delta_j \to \{0\}$ as $j \to \infty$.
\end{enumerate}
\end{definition}

\begin{remark}
This differs from the Dirac delta sequences and families of \cite{Kru} in the sense that no symmetry
$\delta_j(x^{-1}) = \delta_j(x)$ is imposed, matching the relaxation made in the Heisenberg Dirac delta
family of \cite{Khn}. This weaker definition suffices for everything that follows.
\end{remark}

\begin{definition}\label{def:spectral}
A linear map $B:L^1(\Hn)\to C^b(\Omega\times\Cn)$ satisfies the \emph{joint spectral property} if, for every $f\in \mathcal S(\Hn)$, $a=(\lambda,(2k+n)|\lambda|)\in R_k^*$, and $z\in\Cn$,
\begin{enumerate}
    \item[(i)] $B(\mathcal Lf)(a,z)=(2k+n)|\lambda|\,B(f)(a,z)$, and
    \item[(ii)] $B(Tf)(a,z)=i\lambda\,B(f)(a,z)$.
\end{enumerate}
\end{definition}

Note that the zero map trivially satisfies the joint spectral property.

\begin{definition}\label{def:translation}
A linear map $B:L^1(\Hn)\to C^b(\Omega\times\Cn)$ satisfies the \emph{twisted translation property} if
$$
B(L_{(w,s)}f)(a,z)=e^{-i\lambda s}\,e^{-\frac i2\lambda\operatorname{Im}(w\cdot\bar z)}\,B(f)(a,z-w),
$$
for every $(w,s)\in\Hn$, $f\in L^1(\Hn)$, $a=(\lambda,(2k+n)|\lambda|)\in\Omega^*$, and $z\in\Cn$.
\end{definition}

\begin{proposition}\label{prop:spectral-translation}
The Strichartz Fourier transform satisfies
\begin{enumerate}
    \item[\textup{(i)}] the joint spectral property, and
    \item[\textup{(ii)}] the twisted translation property.
\end{enumerate}
\end{proposition}

\begin{proof}
(i) Fix $f\in \mathcal S(\Hn)$ and $a=(\lambda,(2k+n)|\lambda|)\in R_k^*$. The function $e_a(z,t)=e^{i\lambda t}\varphi_{k,\lambda}^{n-1}(z)$ is smooth, with all derivatives bounded, and $\widehat f(a,z)=(f*e_a)(z,0)$. By \Cref{lem:leftright} with $X=T$ and $\widetilde T=T$,
$$
\widehat{Tf}(a,z)=\big((Tf)*e_a\big)(z,0)=\big(f*Te_a\big)(z,0)=i\lambda\,\widehat f(a,z).
$$
Since $e_a$ is radial in $z$, we have $\widetilde{\mathcal L}e_a=\mathcal Le_a=(2k+n)|\lambda|e_a$. Combining this with \Cref{lem:leftright} for $\mathcal L$, we get
\begin{eqnarray*}
\widehat{\mathcal Lf}(a,z) & = & \big((\mathcal Lf)*e_a\big)(z,0)\ =\ \big(f*\widetilde{\mathcal L}e_a\big)(z,0)\\
& = & (2k+n)|\lambda|\,\widehat f(a,z).
\end{eqnarray*}

(ii) Fix $(w,s)\in\Hn$, $f\in L^1(\Hn)$, $a=(\lambda,(2k+n)|\lambda|)\in R_k^*$, and $z\in\Cn$. Write $$
L_{(w,s)}f(u,v)
= f\big((w,s)^{-1}(u,v)\big)
= f\big(u-w,\,v-s-\tfrac12\operatorname{Im}(w\cdot\bar u)\big).
$$
Substituting $(u,v)=(w,s)(u',v')$, we get
\begin{eqnarray*}
\widehat{L_{(w,s)}f}(a,z) & = & \int_{\Hn} f(u',v')\,e_a\big(((w,s)(u',v'))^{-1}(z,0)\big)\,du'\,dv'\\
& = & \int_{\Hn} f(u',v')\,e_a\big((u',v')^{-1}(w,s)^{-1}(z,0)\big)\,du'\,dv'.
\end{eqnarray*}
Since $(w,s)^{-1}(z,0)=\big(z-w,\,-s-\tfrac12\operatorname{Im}(w\cdot\bar z)\big)$, the group law gives
\begin{multline}\label{eq:groupexp}
e_a\big((u',v')^{-1}(z-w,\,-s-\tfrac12\operatorname{Im}(w\cdot\bar z))\big)\\
= e^{-i\lambda s}e^{-\frac i2\lambda\operatorname{Im}(w\cdot\bar z)}\, e^{-i\lambda v'}\,e^{-\frac i2\lambda\operatorname{Im}(u'\cdot\overline{(z-w)})}\varphi_{k,\lambda}^{n-1}\big((z-w)-u'\big).
\end{multline}
Substituting \eqref{eq:groupexp} into the integral and separating the factors independent of $u',v'$,
$$
\widehat{L_{(w,s)}f}(a,z)
= e^{-i\lambda s}e^{-\frac i2\lambda\operatorname{Im}(w\cdot\bar z)}
\int_{\Cn}\left(\int_{\mathbb R} f(u',v')\,e^{-i\lambda v'}\,dv'\right)
e^{-\frac i2\lambda\operatorname{Im}(u'\cdot\overline{(z-w)})}\,\varphi_{k,\lambda}^{n-1}\big((z-w)-u'\big)\,du'.
$$
The inner integral is, by definition, $f^{-\lambda}(u')$, so
\begin{align*}
\widehat{L_{(w,s)}f}(a,z)
&= e^{-i\lambda s}e^{-\frac i2\lambda\operatorname{Im}(w\cdot\bar z)}
\int_{\Cn} f^{-\lambda}(u')\,e^{-\frac i2\lambda\operatorname{Im}(u'\cdot\overline{(z-w)})}\,\varphi_{k,\lambda}^{n-1}\big((z-w)-u'\big)\,du' \\
&= e^{-i\lambda s}e^{-\frac i2\lambda\operatorname{Im}(w\cdot\bar z)}
\int_{\Cn} \varphi_{k,\lambda}^{n-1}\big((z-w)-u'\big)\,f^{-\lambda}(u')\,e^{\frac i2\lambda\operatorname{Im}((z-w)\cdot\bar u')}\,du' \\
&= e^{-i\lambda s}e^{-\frac i2\lambda\operatorname{Im}(w\cdot\bar z)}\, \varphi_{k,\lambda}^{n-1}*_\lambda f^{-\lambda}(z-w) \\
&= e^{-i\lambda s}e^{-\frac i2\lambda\operatorname{Im}(w\cdot\bar z)}\, f^{-\lambda}*_{-\lambda}\varphi_{k,\lambda}^{n-1}(z-w), \ \ (by \eqref{eq:symmetry})\\
&= e^{-i\lambda s}e^{-\frac i2\lambda\operatorname{Im}(w\cdot\bar z)}\,\widehat f(a,z-w),
\end{align*}
which is the stated formula.
\end{proof}

\begin{remark}
Since $\Hn$ is unimodular, $\|L_{(w,s)}f\|_1=\|f\|_1$ for every $(w,s)\in\Hn$, so $L_{(w,s)}f\in L^1(\Hn)$ whenever $f\in L^1(\Hn)$. Hence \Cref{prop:spectral-translation}(ii) is meaningful for every $f\in L^1(\Hn)$.
\end{remark}

\begin{definition}\label{def:strichartz-dirac}
A linear map $B:L^1(\Hn)\to C^b(\Omega\times\Cn)$ satisfies the \emph{Strichartz Dirac delta property} if
$$
\lim_{j\to\infty}B(\delta_j)(a,z)=\varphi_{k,\lambda}^{n-1}(z),\qquad \forall\,z\in\Cn,\ a=(\lambda,(2k+n)|\lambda|)\in\Omega^*,
$$
for some Heisenberg Dirac delta sequence $(\delta_j)$ on $\Hn$.
\end{definition}

\begin{proposition}\label{prop:delta}
The Strichartz Fourier transform satisfies the Strichartz Dirac delta property. In fact, $\widehat{\delta_j}(a,z)\to\varphi_{k,\lambda}^{n-1}(z)$ for every Heisenberg Dirac delta sequence $(\delta_j)$ on $\Hn$.
\end{proposition}
\begin{proof}
Let $(\delta_j)$ be a Heisenberg Dirac delta sequence on $\Hn$, and fix $z\in\Cn$ and $a=(\lambda,(2k+n)|\lambda|)\in R_k^*$. By \eqref{eq:groupexp},
\begin{align*}
\widehat{\delta_j}(a,z)
&=\int_{\Hn}\delta_j(w,s)\,e^{-i\lambda s}e^{-\frac i2\lambda\operatorname{Im}(w\cdot\bar z)}\varphi_{k,\lambda}^{n-1}(z-w)\,dw\,ds\\
&\longrightarrow \varphi_{k,\lambda}^{n-1}(z), \textrm{ \ as \ } \operatorname{supp}\delta_j\to\{(0,0)\},
\end{align*}
since $\int_{\Hn}\delta_j=1$ and the integrand is continuous in $(w,s)$ at $(0,0)$ with value $\varphi_{k,\lambda}^{n-1}(z)$, so the convergence follows from the standard approximate identity argument.
\end{proof}

For $\lambda\in\R^*$, let
$$
\mathcal D_\lambda f := f''+\frac{2n-1}{\rho}f'-\frac{\lambda^2}{4}\rho^2f.
$$
On radial functions $g(z)=g_0(|z|)$, $\mathcal D_\lambda g_0 = -L_\lambda g$, where $L_\lambda$ is the special Hermite operator (\cref{sec:prelims}). To keep computations simpler, we use $\mathcal D_\lambda$ as separate notation, confined to the one-variable setting of this section.

\begin{lemma}\label{lem:selfadjoint}
The operator $\mathcal D_\lambda$ is formally self-adjoint with respect to the measure $\rho^{2n-1}\,d\rho$ on $\R_+$, that is, for every $f,K\in C_c^\infty(\R_+)$,
$$
\int_0^\infty (\mathcal D_\lambda f)(\rho)K(\rho)\rho^{2n-1}\,d\rho=\int_0^\infty f(\rho)(\mathcal D_\lambda K)(\rho)\rho^{2n-1}\,d\rho.
$$
\end{lemma}

\begin{proof}
Integrating by parts (twice on the first term, once on the second),
\begin{align*}
\int_0^\infty f''(\rho)K(\rho)\rho^{2n-1}\,d\rho
&=\int_0^\infty f(\rho)\Big[K''(\rho)\rho^{2n-1}+2(2n-1)K'(\rho)\rho^{2n-2}\\
&\qquad\qquad+(2n-1)(2n-2)K(\rho)\rho^{2n-3}\Big]\,d\rho,
\end{align*}
and
$$
\int_0^\infty \frac{2n-1}{\rho}f'(\rho)K(\rho)\rho^{2n-1}\,d\rho=-(2n-1)\int_0^\infty f(\rho)\Big[K'(\rho)\rho^{2n-2}+(2n-2)K(\rho)\rho^{2n-3}\Big]\,d\rho.
$$
Adding the above two identities, we get
$$
\int_0^\infty\Big(f''+\frac{2n-1}{\rho}f'\Big)(\rho)K(\rho)\rho^{2n-1}\,d\rho=\int_0^\infty f(\rho)\Big(K''+\frac{2n-1}{\rho}K'\Big)(\rho)\rho^{2n-1}\,d\rho.
$$
The multiplication operator $f\mapsto -\tfrac{\lambda^2}4\rho^2f$ is self-adjoint with respect to $\rho^{2n-1}d\rho$, and the lemma follows.
\end{proof}

For $k\in\N$, let $L^1_k:=L^1\big(\R_+,\rho^{2n-1}\,d\rho\big)$.

\begin{definition}\label{def:laguerre}
Fix $k\in\N$. A \emph{bounded} linear map $B_k:L^1_k\to C^b(\Rstar)$ is said to satisfy the \emph{Laguerre differential property} if
$$
B_k\Big(f_0''+\frac{2n-1}{\rho}f_0'-\frac{\lambda^2}{4}\rho^2f_0\Big)(\lambda)=-(2k+n)|\lambda|\,B_k(f_0)(\lambda),
$$
for every $f_0\in C_c^\infty(\R_+)$ and $\lambda\in\Rstar$.
\end{definition}

\begin{proposition}\label{prop:laguerreproperty}
For each $k\in\N$, the map $B_k:L^1_k\to C^b(\Rstar)$ defined by
$$
B_k(f_0)(\lambda) := \int_0^\infty \varphi_{k,\lambda}^{n-1}(\rho)\,f_0(\rho)\,\rho^{2n-1}\,d\rho,\qquad f_0\in L^1_k,\ \lambda\in\Rstar,
$$
satisfies the Laguerre differential property.
\end{proposition}

\begin{proof}
Fix $\lambda\in\Rstar$ and $f_0\in C_c^\infty(\R_+)$. By \eqref{eq:laguerreODE}, and \Cref{lem:selfadjoint} applied with $K=\varphi_{k,\lambda}^{n-1}$,

\begin{align*}
\lefteqn{B_k\Big(f_0''+\tfrac{2n-1}{\rho}f_0'-\tfrac{\lambda^2}{4}\rho^2f_0\Big)(\lambda)}\\
&=\int_0^\infty f_0(\rho)\Big({\varphi_{k,\lambda}^{n-1}}''(\rho)+\tfrac{2n-1}{\rho}{\varphi_{k,\lambda}^{n-1}}'(\rho)-\tfrac{\lambda^2}{4}\rho^2\varphi_{k,\lambda}^{n-1}(\rho)\Big)\rho^{2n-1}\,d\rho\\
&=-(2k+n)|\lambda|\int_0^\infty f_0(\rho)\varphi_{k,\lambda}^{n-1}(\rho)\rho^{2n-1}\,d\rho\\
&=-(2k+n)|\lambda|\,B_k(f_0)(\lambda). \qedhere
\end{align*}
\end{proof}

\section{A weak spectral characterisation}
\label{sec:spectral}

In this section, we obtain a characterisation of the Strichartz Fourier transform on each ray, using the joint spectral property, the twisted translation property, and a normalisation at one explicit test function.

The following result supplies the central-direction piece of the decomposition used in
Step~1 of the proof of \Cref{thm:spectral} below.
\begin{lemma}\label{lem:primitive}
Fix $\lambda\in\Rstar$. If $f\in\mathcal S(\Hn)$ satisfies $f^{-\lambda}=0$, then
$$
f_1(z,t):=e^{i\lambda t}\int_{-\infty}^t e^{-i\lambda s}f(z,s)\,ds
$$
belongs to $\mathcal S(\Hn)$ and satisfies $(T-i\lambda)f_1=f$.
\end{lemma}

\begin{proof}
Fix $z\in\Cn$. Since $f(z,\cdot)\in\mathcal S(\R)$, the function
\begin{equation}\label{eq:F-left}
F(z,t) := \int_{-\infty}^t e^{-i\lambda s}f(z,s)\,ds
\end{equation}
is smooth in $t$, with $\partial_tF(z,t)=e^{-i\lambda t}f(z,t)$. Since $f^{-\lambda}(z)=\int_\R e^{-i\lambda s}f(z,s)\,ds=0$,
\begin{equation}\label{eq:F-right}
F(z,t) = -\int_t^\infty e^{-i\lambda s}f(z,s)\,ds.
\end{equation}
Rapid decrease of $f(z,\cdot)$, applied to \eqref{eq:F-left} as $t\to-\infty$ and to \eqref{eq:F-right} as $t\to+\infty$, gives rapid decay of $F(z,\cdot)$ and of every $t$-derivative of $F(z,\cdot)$. Differentiating under the integral sign in $z$ gives the same bounds uniformly on compact $z$-sets, and similarly for $z$-derivatives. Thus $F\in\mathcal S(\Hn)$, and hence so is
$$
f_1(z,t) = e^{i\lambda t}F(z,t).
$$
Finally,
\begin{align*}
(T-i\lambda)f_1 &= \partial_t\big(e^{i\lambda t}F\big)-i\lambda e^{i\lambda t}F
= e^{i\lambda t}\big(i\lambda F+\partial_tF\big)-i\lambda e^{i\lambda t}F\\
&= e^{i\lambda t}\partial_tF
= e^{i\lambda t}e^{-i\lambda t}f = f. \qedhere
\end{align*}
\end{proof}

For $\lambda\in\Rstar$, let
$$
\chi_\lambda(t):=e^{\lambda^2/4\pi}e^{-\pi t^2},\qquad t\in\R.
$$
Then $\chi_\lambda\in\mathcal S(\R)$ and $\int_\R\chi_\lambda(t)e^{-i\lambda t}\,dt=1$.

The following result is our first characterisation of the Strichartz Fourier transform.

\begin{theorem}\label{thm:spectral}
Fix $k\in\N$. Let $B:L^1(\Hn)\to C^b(\Omega\times\Cn)$ be a linear, bounded map which satisfies
\begin{enumerate}
    \item[\textup{(i)}] the joint spectral property,
    \item[\textup{(ii)}] the twisted translation property, and
    \item[\textup{(iii)}] $B\big(\varphi_{k,\lambda}^{n-1}\otimes\chi_\lambda\big)(a,z)=(2\pi)^n|\lambda|^{-n}\varphi_{k,\lambda}^{n-1}(z)$, for every $a=(\lambda,(2k+n)|\lambda|)\in R_k^*$ and $z\in\Cn$.
\end{enumerate}
Then $B(f)(a,z)=\widehat f(a,z)$ for every $f\in L^1(\Hn)$, every $a\in R_k^*$, and every $z\in\Cn$.
\end{theorem}

\begin{proof}
Fix $a=(\lambda,(2k+n)|\lambda|)\in R_k^*$, and write $$\Phi_a(f):=\widehat f(a,\cdot)=f^{-\lambda}*_{-\lambda}\varphi_{k,\lambda}^{n-1}$$ for $f\in\mathcal S(\Hn)$.

Fix $g\in\mathcal S(\Cn)$ and set $f:=g\otimes\chi_\lambda\in\mathcal S(\Hn)$. Since $f^{-\lambda}=g\cdot\int_\R\chi_\lambda(t)e^{-i\lambda t}\,dt=g$,
$$
\Phi_a(f)=g*_{-\lambda}\varphi_{k,\lambda}^{n-1}=(2\pi)^n|\lambda|^{-n}P_k(g).
$$
As $g$ ranges over $\mathcal S(\Cn)$, and $P_k$ carries $\mathcal S(\Cn)$ into $\mathcal S(\Cn)\cap\mathcal E_k^\lambda$ (\Cref{sec:prelims}). Note that every $h\in\mathcal S(\Cn)\cap\mathcal E_k^\lambda$ arises this way by taking $g=h$, as  $P_k(h)=h$. Thus we obtain
\begin{equation}\label{eq:step1range}
\Phi_a\big(\mathcal S(\Hn)\big)=(2\pi)^n|\lambda|^{-n}\big(\mathcal S(\Cn)\cap\mathcal E_k^\lambda\big),
\end{equation}
a subspace dense in $\mathcal E_k^\lambda$, since $\mathcal S(\Cn)$ is dense in $L^2(\Cn)$ and $P_k$ is continuous.

Next, suppose $f\in\mathcal S(\Hn)$ and $\Phi_a(f)=0$. This gives \ $g:=f^{-\lambda}\perp\mathcal E_k^\lambda$. The resolvent $(L_{-\lambda}-(2k+n)|\lambda|)^{-1}$ is bounded on $(\mathcal E_k^\lambda)^\perp$ and preserves $\mathcal S(\Cn)$ (\Cref{sec:prelims}). Taking
\begin{align*}
g_1&:=\big(L_{-\lambda}-(2k+n)|\lambda|\big)^{-1}g\in\mathcal S(\Cn),
\end{align*}
we get $(L_{-\lambda}-(2k+n)|\lambda|)g_1=g$.\\
Let $f_2:=g_1\otimes\chi_\lambda\in\mathcal S(\Hn)$. Since $(\mathcal Lh)^\mu=L_\mu(h^\mu)$ for every $\mu$,
$$
\big((\mathcal L-(2k+n)|\lambda|)f_2\big)^{-\lambda}=\big(L_{-\lambda}-(2k+n)|\lambda|\big)g_1=g=f^{-\lambda}.
$$
So $f-(\mathcal L-(2k+n)|\lambda|)f_2\in\mathcal S(\Hn)$ has vanishing $(-\lambda)$-th central Fourier coefficient. By \Cref{lem:primitive},
\begin{align*}
f-\big(\mathcal L-(2k+n)|\lambda|\big)f_2 &= (T-i\lambda)f_1,
\end{align*}
and hence
\begin{align*}
f &= \big(\mathcal L-(2k+n)|\lambda|\big)f_2+(T-i\lambda)f_1.
\end{align*}

Conversely, by \Cref{prop:spectral-translation}(i) applied to $\widehat{(\cdot)}$ itself, $\Phi_a$ vanishes on $(\mathcal L-(2k+n)|\lambda|)f_2$ and on $(T-i\lambda)f_1$, for every $f_1,f_2\in\mathcal S(\Hn)$. So $\ker\big(\Phi_a|_{\mathcal S(\Hn)}\big)$ equals the set of such sums.

The Strichartz Fourier transform satisfies the joint spectral property (\Cref{prop:spectral-translation}(i)), and $B$ satisfies it by hypothesis. So both $f\mapsto\widehat f(a,\cdot)$ and $f\mapsto B(f)(a,\cdot)$ vanish on $\ker(\Phi_a|_{\mathcal S(\Hn)})$, as just characterised. Consequently, whenever $f,f'\in\mathcal S(\Hn)$ satisfy $\Phi_a(f)=\Phi_a(f')$, we have
\begin{equation}\label{eq:welldef}
B(f)(a,\cdot)=B(f')(a,\cdot).
\end{equation}

 Write $\phi:=\varphi_{k,\lambda}^{n-1}$. For $u\in\mathcal S(\Cn)\cap\mathcal E_k^\lambda$, let $f_u:=u\otimes\chi_\lambda\in\mathcal S(\Hn)$, so that, by \eqref{eq:step1range}, we have $\Phi_a(f_u)=(2\pi)^n|\lambda|^{-n}u$. Let
$$
F_u:=f_u*(\phi\otimes\chi_\lambda)\in\mathcal S(\Hn),
$$
 where $*$ denotes convolution on $\Hn$. Writing out the convolution using the group law, as in the proof of \Cref{prop:spectral-translation}(ii),
$$
F_u(z,t)=\int_{\Cn} u(w)\,\phi(z-w)\,(\chi_\lambda*\chi_\lambda)\Big(t-\tfrac12\operatorname{Im}(w\cdot\bar z)\Big)\,dw,
$$
where $\chi_\lambda*\chi_\lambda$ denotes ordinary convolution on $\R$. Since
$$
\int_\R(\chi_\lambda*\chi_\lambda)(t)e^{-i\lambda t}\,dt=\Big(\int_\R\chi_\lambda(t)e^{-i\lambda t}\,dt\Big)^2=1.
$$
Taking the $(-\lambda)$-th central Fourier coefficient, and expanding the kernel of the twisted convolution it reduces to, 
\begin{equation}\label{eq:step3closing}
\begin{aligned}
F_u^{-\lambda}(z) & = \int_{\Cn}u(w)\,\phi(z-w)\,e^{-\frac i2\lambda\operatorname{Im}(w \cdot \bar z)}\,dw \ =\ \phi*_\lambda u(z)\\
& = u*_{-\lambda}\phi(z)\ =\ (2\pi)^n|\lambda|^{-n}u(z),
\end{aligned}
\end{equation}
using \eqref{eq:symmetry} and \eqref{eq:special-case}, since $u\in\mathcal E_k^\lambda$. Hence
$$
\Phi_a(F_u)=F_u^{-\lambda}*_{-\lambda}\phi=(2\pi)^n|\lambda|^{-n}u*_{-\lambda}\phi=(2\pi)^{2n}|\lambda|^{-2n}u.
$$
 Since $F_u=\int_{\Hn}f_u(y)\,L_y(\phi\otimes\chi_\lambda)\,dy$, $B$ is linear and bounded, and $y\mapsto L_y(\phi\otimes\chi_\lambda)$ is continuous from $\Hn$ into $L^1(\Hn)$, we get
$$
B(F_u)(a,z)=\int_{\Hn}f_u(w,s)\,B\big(L_{(w,s)}(\phi\otimes\chi_\lambda)\big)(a,z)\,dw\,ds.
$$
By the twisted translation property,
\begin{align*}
B\big(L_{(w,s)}(\phi\otimes\chi_\lambda)\big)(a,z)
&=e^{-i\lambda s}e^{-\frac i2\lambda\operatorname{Im}(w\cdot\bar z)}B(\phi\otimes\chi_\lambda)(a,z-w).
\end{align*}
By hypothesis (iii), $B(\phi\otimes\chi_\lambda)(a,z-w)=(2\pi)^n|\lambda|^{-n}\phi(z-w)$. Substituting,
\begin{eqnarray*}
B(F_u)(a,z) & = & \Big(\int_\R\chi_\lambda(s)e^{-i\lambda s}\,ds\Big)(2\pi)^n|\lambda|^{-n}\int_{\Cn}u(w)\,e^{-\frac i2\lambda\operatorname{Im}(w\cdot\bar z)}\phi(z-w)\,dw\\
& = & (2\pi)^n|\lambda|^{-n}\big[\phi*_\lambda u\big](z).
\end{eqnarray*}
As in \eqref{eq:step3closing}, $\phi*_\lambda u=u*_{-\lambda}\phi=(2\pi)^n|\lambda|^{-n}u$, so
\begin{equation}\label{eq:step4result}
B(F_u)(a,z)=(2\pi)^{2n}|\lambda|^{-2n}u(z).
\end{equation}

 Comparing \eqref{eq:step1range}, \eqref{eq:step3closing}, and \eqref{eq:step4result},
 \begin{align*}
B(F_u)(a,\cdot) &= \Phi_a(F_u)\\
&= (2\pi)^{2n}|\lambda|^{-2n}u\\
&= (2\pi)^n|\lambda|^{-n}\Phi_a(f_u)\\
&= \Phi_a\big((2\pi)^n|\lambda|^{-n}f_u\big)
\end{align*}
Thus
 \begin{align*}
\Phi_a\Big(F_u-(2\pi)^n|\lambda|^{-n}f_u\Big)=0.
\end{align*}
By \eqref{eq:welldef}, $B\big(F_u-(2\pi)^n|\lambda|^{-n}f_u\big)(a,\cdot)=0$, that is, $B(F_u)(a,\cdot)=(2\pi)^n|\lambda|^{-n}B(f_u)(a,\cdot)$. Combined with \eqref{eq:step4result},
\begin{eqnarray*}
(2\pi)^n|\lambda|^{-n}B(f_u)(a,z) & = & (2\pi)^{2n}|\lambda|^{-2n}u(z),\\
\text{so}\qquad B(f_u)(a,z) & = & (2\pi)^n|\lambda|^{-n}u(z)=\Phi_a(f_u)(z)=\widehat{f_u}(a,z).
\end{eqnarray*}
By \eqref{eq:step1range}, 
$$\Phi_a(\mathcal S(\Hn))=\{(2\pi)^n|\lambda|^{-n}u:u\in\mathcal S(\Cn)\cap\mathcal E_k^\lambda\},$$ so every $f\in\mathcal S(\Hn)$ has $\Phi_a(f)=\Phi_a(f_u)$ for some such $u$. Then, by \eqref{eq:welldef} again, $$B(f)(a,\cdot)=B(f_u)(a,\cdot)=\widehat{f_u}(a,\cdot)=\Phi_a(f)=\widehat f(a,\cdot).$$ Thus $B(f)(a,z)=\widehat f(a,z)$ for every $f\in\mathcal S(\Hn)$. Since $\mathcal S(\Hn)$ is dense in $L^1(\Hn)$, and both $B$ and $f\mapsto\widehat f(a,\cdot)$ are continuous from $L^1(\Hn)$ to $C^b(\Cn)$ (\cite{Th23}), this extends to every $f\in L^1(\Hn)$.
\end{proof}

\section{A characterisation via translation}
\label{sec:translation}

This section discusses a translation-based characterisation of the Strichartz Fourier transform.

\begin{theorem}\label{thm:translation}
Let $B:L^1(\Hn)\to C^b(\Omega\times\Cn)$ be a bounded linear map which satisfies the twisted translation property and the Strichartz Dirac delta property. Then $B(f)(a,z)=\widehat f(a,z)$ for every $f\in L^1(\Hn)$, every $a\in\Omega^*$, and every $z\in\Cn$.
\end{theorem}

\begin{proof}
Let $(\delta_j)$ be a Heisenberg Dirac delta sequence on $\Hn$ such that $B$ satisfies the Strichartz Dirac delta property. Then for every $f\in L^1(\Hn)$, we have
 $f*\delta_j\to f$ in $L^1(\Hn)$ as $j\to\infty$ (\cite{Folland}).
 
Fix $f\in L^1(\Hn)$, $a=(\lambda,(2k+n)|\lambda|)\in R_k^*$, and $z\in\Cn$. Since $B$ is a bounded linear map, and $y\mapsto L_yf$ is continuous from $\Hn$ into $L^1(\Hn)$,
$$
B(f*\delta_j)(a,z)=\int_{\Hn} f(w,s)\,B\big(L_{(w,s)}\delta_j\big)(a,z)\,dw\,ds.
$$
Since $B$ satisfies the twisted translation property,
$$
B\big(L_{(w,s)}\delta_j\big)(a,z)=e^{-i\lambda s}e^{-\frac i2\lambda\operatorname{Im}(w\cdot\bar z)}B(\delta_j)(a,z-w),
$$
we get
\begin{align}
\label{eq:bfstard}
B(f*\delta_j)(a,z)=\int_{\Hn} f(w,s)\,e^{-i\lambda s}e^{-\frac i2\lambda\operatorname{Im}(w\cdot\bar z)}\,B(\delta_j)(a,z-w)\,dw\,ds.
\end{align}
As $j\to\infty$, the LHS of \eqref{eq:bfstard} gives
\begin{align}
\label{eq:bfazlim}
B(f*\delta_j)(a,z)\longrightarrow B(f)(a,z),
\end{align}
since $f*\delta_j\to f$ in $L^1(\Hn)$ and $B$ is bounded.\\ 
For each $j$, the RHS of \eqref{eq:bfstard} gives
$$
|B(\delta_j)(a,z-w)|\le\|B\|_{\mathrm{op}}\|\delta_j\|_1=\|B\|_{\mathrm{op}}, \qquad \forall\, j,w,s.
$$
Since $B$ satisfies the Strichartz Dirac delta property,we get $$B(\delta_j)(a,z-w)\to\varphi_{k,\lambda}^{n-1}(z-w) \textrm{ \ pointwise in \ } (w,s).$$ Since $f\in L^1(\Hn)$, the dominated convergence theorem gives
$$
B(f*\delta_j)(a,z)\longrightarrow \int_{\Hn} f(w,s)\,e^{-i\lambda s}e^{-\frac i2\lambda\operatorname{Im}(w\cdot\bar z)}\,\varphi_{k,\lambda}^{n-1}(z-w)\,dw\,ds.
$$
By \eqref{eq:groupexp}, the integrand on the right is $f(w,s)\,e_a\big((w,s)^{-1}(z,0)\big)$, so this limit equals $\widehat f(a,z)$.\\
Combining this with \eqref{eq:bfazlim} gives
$$
B(f)(a,z)=\widehat f(a,z).
$$
\end{proof}

The following remark compares some of the results in this paper among themselves.
\begin{remark}
\label{rem:eqlnce}
In the proof of \Cref{thm:translation}, the twisted translation property determines $B(\delta_j)(a,z-w)$, for every $w\in\Cn$, from $B(L_{(w,s)}\delta_j)(a,z)$. When combined with the Strichartz Dirac delta property, this already determines $B$ on its own, without any spectral hypothesis. Thus, \Cref{thm:translation} is not a special case of \Cref{thm:spectral}. On the other hand, \Cref{thm:spectral} replaces the Strichartz Dirac delta property by the joint spectral property together with a single algebraic identity at one explicit test function. (Note that the Strichartz Dirac delta property is a limiting condition along an entire approximate identity.)
\end{remark}

\section{A Laguerre differential characterisation}
\label{sec:laguerre}

Now we obtain a characterisation of the Strichartz Fourier transform in terms of the Laguerre differential property, and then deduce a characterisation of the Hankel transform.

First we obtain that the target map is an integral transform, with an appropriately bounded kernel.

\begin{proposition}\label{prop:kernelfromboundedness}
Let $B_k:L^1_k\to C^b(\Rstar)$ be a bounded linear map satisfying the Laguerre differential property. Then for every $\lambda\in\Rstar$ there is a function $K_k(\lambda,\cdot)\in C^\infty(0,\infty)$ such that
$$
B_k(f_0)(\lambda)=\int_0^\infty K_k(\lambda,\rho)\,f_0(\rho)\,\rho^{2n-1}\,d\rho,\qquad \forall\,f_0\in L^1_k,
$$
$K_k(\lambda,\cdot)$ satisfies \eqref{eq:laguerreODE} classically on $(0,\infty)$, and
$$
\sup_{\rho>0}|K_k(\lambda,\rho)|\le\|B_k\|_{\mathrm{op}}.
$$
\end{proposition}

\begin{proof} Fix $\lambda\in\Rstar$. The map $f_0\mapsto B_k(f_0)(\lambda)$ is a bounded linear functional on $L^1_k=L^1(\R_+,\rho^{2n-1}d\rho)$ satisfying
$$
|B_k(f_0)(\lambda)|\le\|B_k(f_0)\|_{C^b(\Rstar)}\le\|B_k\|_{\mathrm{op}}\|f_0\|_{L^1_k}.
$$
Since $(\R_+,\rho^{2n-1}d\rho)$ is $\sigma$-finite, by $L^1$-$L^\infty$ duality (\cite{Rud}), there exists a unique $K_k(\lambda,\cdot)\in L^\infty(\R_+,\rho^{2n-1}d\rho)$ with
$$
B_k(f_0)(\lambda)=\int_0^\infty K_k(\lambda,\rho)f_0(\rho)\rho^{2n-1}\,d\rho,\qquad \forall\,f_0\in L^1_k,
$$
\begin{equation}\label{eq:aeubound}
\operatorname{ess\,sup}_{\rho>0}|K_k(\lambda,\rho)|\le\|B_k\|_{\mathrm{op}}.\end{equation}
 Let $u:=K_k(\lambda,\cdot)\in L^\infty(0,\infty)\subset L^1_{\mathrm{loc}}(0,\infty)$.
 
Since
$$
\rho^{2n-1}\,\mathcal D_\lambda f_0=(\rho^{2n-1}f_0')'-\frac{\lambda^2}{4}\rho^{2n+1}f_0,\qquad f_0\in C_c^\infty(\R_+),
$$
the Laguerre differential property, written against this form, gives
\begin{equation}\label{eq:weakform}
\int_0^\infty u\,(\rho^{2n-1}f_0')'\,d\rho-\frac{\lambda^2}{4}\int_0^\infty u\,\rho^{2n+1}f_0\,d\rho
=-(2k+n)|\lambda|\int_0^\infty u\,\rho^{2n-1}f_0\,d\rho
\end{equation}
for every $f_0\in C_c^\infty(\R_+)$. Let $\psi:=\rho^{2n-1}f_0'\in C_c^\infty(\R_+)$. The definitions of the distributional derivative (\cite{Rud}) of $u$ and of $\rho^{2n-1}u'$ give
$$
\int_0^\infty u\,\psi'\,d\rho=-\langle u',\psi\rangle=-\langle \rho^{2n-1}u',f_0'\rangle=\langle(\rho^{2n-1}u')',f_0\rangle,
$$
with no regularity of $u$ assumed at any step. Substituting into \eqref{eq:weakform} gives the following identity of distributions on $(0,\infty)$,
\begin{equation}\label{eq:distform}
(\rho^{2n-1}u')'=g,\qquad g:=\frac{\lambda^2}{4}\rho^{2n+1}u-(2k+n)|\lambda|\rho^{2n-1}u.
\end{equation}
Since $u\in L^\infty(0,\infty)$ and $\rho^{2n+1},\rho^{2n-1}$ are smooth, we get $g\in L^\infty_{\mathrm{loc}}(0,\infty)\subset L^1_{\mathrm{loc}}(0,\infty)$.

We show, by induction on $m$, that
$$
u\in C^{2m+1}[\varepsilon,M],\qquad \forall\,m\ge0,\ \ \forall\,0<\varepsilon<M<\infty.
$$
\emph{Base case ($m=0$).} Since $g\in L^1_{\mathrm{loc}}[\varepsilon,M]$, \eqref{eq:distform} shows that $\rho^{2n-1}u'$ agrees a.e.\ with a continuous function on $[\varepsilon,M]$. As $\rho^{2n-1}$ is smooth and nonvanishing on $[\varepsilon,M]$, the map
$$
u'=\rho^{1-2n}\cdot(\rho^{2n-1}u')
$$
is continuous there, so $u\in C^1[\varepsilon,M]$.

\emph{Inductive step.} Suppose $u\in C^{2m+1}[\varepsilon,M]$. Then $g\in C^{2m+1}[\varepsilon,M]$. From \eqref{eq:distform} we get $\rho^{2n-1}u'\in C^{2m+2}[\varepsilon,M]$. Dividing by the smooth, nonvanishing function $\rho^{2n-1}$, we get
$$
u'\in C^{2m+2}[\varepsilon,M],\qquad\text{hence}\qquad u\in C^{2m+3}[\varepsilon,M].
$$

Letting $\varepsilon\to0^+$ and $M\to\infty$, we get that $
u\in C^\infty(0,\infty),$
Thus \eqref{eq:distform} (equivalently \eqref{eq:laguerreODE} with eigenvalue $(2k+n)|\lambda|$) now holds on $(0,\infty)$ in the classical sense. Thus the map $u$ has a continuous representative on $(0,\infty)$.

By \eqref{eq:aeubound},
$$
|K_k(\lambda,\rho)| = |u(\rho)|\le\|B_k\|_{\mathrm{op}}\qquad\text{for a.e.\ }\rho>0.
$$
A continuous function that is bounded a.e.\ by $\|B_k\|_{\mathrm{op}}$ is bounded by $\|B_k\|_{\mathrm{op}}$ everywhere on $(0,\infty)$. Hence
$$
\sup_{\rho>0}|K_k(\lambda,\rho)|\le\|B_k\|_{\mathrm{op}}.
$$
\end{proof}

\begin{theorem}\label{thm:laguerre}
Fix $k\in\N$. Let $B_k:L^1_k\to C^b(\Rstar)$ be a bounded linear map satisfying the Laguerre differential property. Then there is a function $c:\Rstar\to\mathbb C$ such that
$$
B_k(f_0)(\lambda)=c(\lambda)\int_0^\infty \varphi_{k,\lambda}^{n-1}(\rho)\,f_0(\rho)\,\rho^{2n-1}\,d\rho,\qquad \forall\,f_0\in L^1_k.
$$
\end{theorem}

\begin{proof}
Fix $\lambda\in\Rstar$. By \Cref{prop:kernelfromboundedness}, $B_k$ is represented by a kernel
$$
K_k(\lambda,\cdot)\in C^\infty(0,\infty)\subset C^2(0,\infty),
$$
satisfying \eqref{eq:laguerreODE} classically on $(0,\infty)$ and bounded near $\rho=0$. This is a second order linear ordinary differential equation with a regular singular point at $\rho=0$, so its solution space is two dimensional, with indicial equation at $\rho=0$
$$
r(r-1)+(2n-1)r=0,\qquad\text{roots } r=0,\ 2-2n.
$$
The root $0$ gives the solution $\varphi_{k,\lambda}^{n-1}$, regular at the origin.

The second, linearly independent, solution $\varphi_2$ is unbounded near $\rho=0$ for every $n\ge1$:
$$
\varphi_2(\rho)\sim
\begin{cases}
\log\rho, & n=1,\\[4pt]
\rho^{2-2n}, & n\ge2.
\end{cases}
$$
By reduction of order,
$$
\varphi_2(\rho) = \varphi(\rho)\int \rho^{1-2n}\varphi(\rho)^{-2}\,d\rho,
$$
and since $\varphi(\rho)^{-2}$ is an even analytic function of $\rho$ near the origin with $\varphi(0)\neq0$,
$$
\rho^{1-2n}\varphi(\rho)^{-2} = a_0\rho^{1-2n} + a_1\rho^{3-2n}+\cdots.
$$
Integrating term by term,
$$
\int \rho^{1-2n}\varphi(\rho)^{-2}\,d\rho =
\begin{cases}
a_0\log\rho + O(1), & n=1,\\[4pt]
\dfrac{a_0}{2-2n}\rho^{2-2n} + \cdots + a_{n-1}\log\rho + O(1), & n\ge 2,
\end{cases}
$$
the $\log\rho$ term for $n\ge2$ arising from the exponent $-1$ recurring at the $(n-1)$-th term of the expansion. In either case, $\varphi_2$ is unbounded near $\rho=0$.

Boundedness of $K_k(\lambda,\cdot)$ near the origin therefore forces
$$
K_k(\lambda,\rho)=c(\lambda)\varphi_{k,\lambda}^{n-1}(\rho),\qquad\text{for some } c(\lambda)\in\mathbb C \text{ not depending on } \rho.
$$
This proves the theorem.
\end{proof}

\begin{corollary}\label{cor:laguerre}
In addition to the hypotheses of \Cref{thm:laguerre}, suppose that for every $\lambda\in\Rstar$ there is a function $f_\lambda\in C_c^\infty(\R_+)$ such that
$$
B_k(f_\lambda)(\lambda)=\int_0^\infty\varphi_{k,\lambda}^{n-1}(\rho)f_\lambda(\rho)\rho^{2n-1}\,d\rho\ne0.
$$
Then
$$
B_k(f_0)(\lambda)=\int_0^\infty\varphi_{k,\lambda}^{n-1}(\rho)\,f_0(\rho)\,\rho^{2n-1}\,d\rho,\qquad\forall\,f_0\in L^1_k,\ \lambda\in\Rstar.
$$
\end{corollary}

\begin{proof}
By \Cref{thm:laguerre},
$$
B_k(f_\lambda)(\lambda)=c(\lambda)\int_0^\infty\varphi_{k,\lambda}^{n-1}(\rho)f_\lambda(\rho)\rho^{2n-1}\,d\rho.
$$
Comparing with the hypothesis gives $c(\lambda)=1$ for every $\lambda\in\Rstar$, since the integral is nonzero.
\end{proof}

\begin{remark}
In the above corollary, the test function $f_\lambda$ must be allowed to depend on $\lambda$. For a single fixed $f_*$, the function
$$
\lambda\longmapsto\int_0^\infty\varphi_{k,\lambda}^{n-1}(\rho)\,f_*(\rho)\,\rho^{2n-1}\,d\rho
$$
need not be free of zeros. Therefore, a normalisation at one test function cannot hold for every $\lambda$ in general. In \Cref{thm:kru5}, a normalisation at a single point suffices, as the transform variable enters the kernel only through the product $y\rho$. However, in the present context, $\lambda$ appears directly in the coefficients of the differential operator in the Laguerre differential property, as a fixed parameter with no rescaling available. Hence, the self-adjointness computation of \Cref{lem:selfadjoint} is simpler than the analogous one for the Bessel operator.
\end{remark}

\begin{remark}
\label{rem:hankel}
We record the degenerate limit of \Cref{thm:laguerre}, in which it reduces to the classical Hankel transform characterisation of \Cref{thm:kru5}. As $k\to\infty$ and $\lambda\to0$ with $(2k+n)|\lambda|\to\tau>0$ fixed, equation \eqref{eq:laguerreODE} degenerates to
$$
-\varphi''-\frac{2n-1}{\rho}\varphi'=\tau\varphi.
$$
A solution regular at the origin is given by
$$
\rho\longmapsto\frac{J_{n-1}(\sqrt\tau\rho)}{(\sqrt\tau\rho)^{n-1}}.
$$
This is the Hankel transform of order $n-1$, consistent with the fact, noted in \cite{Th23}, that the Strichartz Fourier transform reduces to a Hankel-type transform on the limiting ray $R_\infty$ of the Heisenberg fan.
\end{remark}

\section{Concluding remarks}
Our three main results identify the Strichartz Fourier transform from different, non-nested combinations of structural properties: \Cref{thm:spectral} and \Cref{thm:translation} both recover $B=\widehat{\ \cdot\ }$ on each ray, but exchange a single-point normalisation for the joint spectral property of $(\mathcal L,T)$ in opposite directions (\Cref{rem:eqlnce}). \Cref{thm:laguerre} gives a characterisation in terms of the Laguerre differential property. From this, we deduced in \Cref{rem:hankel}, that \Cref{thm:laguerre} degenerates to the classical Hankel transform characterisation of Theorem \ref{thm:kru5} on the limiting ray $R_\infty$. Together, these results give a structural account of the Strichartz Fourier transform across the whole Heisenberg fan $\Omega$ - ray by ray, and at the level of its radial kernel, built entirely from properties intrinsic to $\Hn$. Although both \cite{Kru} and the results above concern scalar-valued transforms, the methods used here are quite distinct, built from properties intrinsic to $\Hn$ rather than adapted directly from \cite{Kru}'s proofs.\\

\end{document}